\documentclass[12pt,a4paper,oneside]{article}
\usepackage{amsfonts,amsmath,amsthm,amssymb,amscd}
\usepackage{fullpage,epsfig}

\usepackage{verbatim}
\newcommand{\beq}{\begin{equation}}
\newcommand{\eeq}{\end{equation}}
\newcommand{\bea}{\begin{eqnarray}}
\newcommand{\eea}{\end{eqnarray}}
\newcommand{\ba}{\begin{array}}
\newcommand{\ea}{\end{array}}

\usepackage{amsmath,amsfonts,amssymb,amsthm}
\usepackage{xspace}
\usepackage{color}

\newtheorem{example}{Example}

\begin{document}

\title{Spatiotemporal Orthogonal Polynomial Approximation for Partial Differential Equations}
\author{ Samir Kumar Bhowmik$^{1, 3}$ and Sharanjeet Dhawan$^2$\\
$1.$ Department of Mathematics and Statistics,  College of Science\\ Al Imam Mohammad Ibn Saud Islamic University(IMSIU),  P.O. BOX 90950\\ 11623 Riyadh, Kingdom of Saudi Arabia,  bhowmiksk@gmail.com\\
$2.$ Department of Mathematics, D.A.V.University, Jalandhar-144012,India\\dhawan311@gmail.com
\\
$3.$ Department of Mathematics, University of Dhaka, Dhaka 1000, Bangladesh\\  bhowmiksk@du.ac.bd
}
\maketitle
{\bf ABSTRACT}

Starting with some fundamental concepts, in this article we present the essential aspects of spectral methods and
their applications to the numerical solution of Partial Differential Equations (PDEs).
We start by using Lagrange and Techbychef orthogonal polynomials  for spatiotemporal approximation of PDEs as a weighted
sum of polynomials. We use collocation at some clustered
grid points to generate a system of equations to approximate the weights for the polynomials.
We finish the study by demonstrating approximate solutions of some  PDEs in one space dimension.

\section{Introduction}
In almost all the areas of science and engineering, a variety of problems are being modeled with the
help of differential equations(ODEs/PDEs). Due to the availability of high speed computers and their
ability to handle large calculations, simulation of these models is feasible.
When we talk about recent advances in computational schemes, we come across: Finite Difference Method,
Finite Element Method, Finite Volume Method, Spectral Method and so on. In Spectral methods,
basis functions are infinitely differentiable and non vanishing on the entire domain (having global support).
On the other hand, bases functions used in finite difference or finite element methods have small compact
support and poor continuity properties. Conclusively, spectral method have good accuracy.
Thus, nowadays, Spectral methods are extensively being used for using higher
order polynomials approximation. Moreover, Spectral Methods achieve a
greater precision with a smaller number of points than Finite Difference methods.

In general, Spectral methods \cite{Trefethen} are a class of spatial discretizations for differential equations.
They can be categorised as Galerkin, tau and collocation(or pseudospectral) spectral methods.
Galerkin and tau work with the coefficients of a global
expansion whereas pseudospectral work with the values at collocation points
and helps in handling the nonlinear terms easily. There are two key components for the formulation of Spectral methods:
\begin{itemize}
  \item Trial function, which is also called the expansion or approximation functions.
  \item Test function, which is also known as weight functions.
\end{itemize}

The trial functions which are linear combination of suitable trial basis functions,
are used to provide the approximate representation of the solution. The test functions
are used to ensure that the differential equations and perhaps some boundary
conditions are satisfied as closely as possible by the truncated series expansion.
Spectral methods are distinguished not only by the fundamental type of the method i.e.
Galerkin, collocation, Galarkin with numerical integration or tau but also by the
particular choice of the trial functions. The most frequently used trial functions are
Trigonometric polynomials, Chebyshev polynomials and Legendre polynomials.
In the present work, we are interested to use orthogonal polynomials $\phi_{j}(x)$
as trial function in the approximation
$
 u_{h}=\sum_{j=0}^{N} a_{j} \phi_{j}.
$
Here, if the function $u$ belongs to $C^{\infty}$ class, the produced error of
approximation as $N$ tends to infinity, approaches zero with exponential rate.
This phenomenon is usually referred to as ``spectral accuracy''.
The accuracy of derivatives obtained by direct, term by term differentiation of
such truncated expansion naturally deteriorates \cite{Atkinson}, but for low-order
derivatives and sufficiently high-order truncations this deterioration is negligible.
So, if solution function and coefficient functions are analytic on the given interval,
the spectral methods will be very efficient and suitable.

 There are evidences of using global polynomials for spatial approximation followed by various one / multistep methods for temporal approximations~\cite{Trefethen}. Because of the stability issues of the temporal scheme one needs small time steps to generate accurate solutions. Thus these solvers are also computationally expensive.

 In this short communication we apply an orthogonal polynomial approximation for both space and time to speed up the computations. This technique  outperforms any existing methods of lines for time integration coupled with any scheme for spatial approximation. Some test problem are studied to check the efficiency of the proposed scheme.

A brief outline of this paper is as follows. In Section~\ref{preli}, we formulate
the spatial and temporal approximation applied to the solution of \ref{appr:f}.
In continuation to Section~\ref{preli}, Section~\ref{compt} is dedicated to computation of
unknown parameters arising in the proposed numerical scheme. Numerical experiments
are reported in Section~\ref{exprmt} for some test problems.
Section \ref{concl} gives a brief summary of the work done and further aspects in our future study.

\section{Polynomial approximation in space and time}\label{preli}
Most popular linear parabolic and hyperbolic PDES can be generalized by
\begin{equation}\label{appr:f}
  \alpha \frac{\partial^2 u}{\partial t^2}+\beta\frac{\partial u}{\partial t} =  \mathbb{L}u + f(x, t)
\end{equation}
where
\[
 \mathbb{L}u = \varepsilon\frac{\partial^2 u}{\partial x^2} + \mu\frac{\partial u}{\partial x}+ \nu u(x, t),
\]
for appropriate choices of real parameters $\alpha$, $\beta$, $ \varepsilon$, $\mu$ and $\nu$. Here we consider the following initial and
 boundary conditions:\\
  \[
  u(\pm 1) = 0, ~~  u(x, 0) = g_0(x),\  \text{and} \  u_t(x,0)=g_1(x).
  \]
We approximate the solution of (\ref{appr:f}) using spectral collocation method. We approximate $u(x, t)$ of (\ref{appr:f}) by a truncated series
\begin{equation}\label{appr1:f}
u_{h} (x,t)=\sum_{i=0}^N\sum_{j=0}^{M} C_{i,j} \psi_i(t)\phi_{j} (x)
\end{equation}
so that $u_h$ satisfies (\ref{appr:f}). Where $C_{i,j}$ are unknown
parameters, $\phi_{j} (x)$ are orthogonal polynomials of suitable type
in $[-1,\ 1]$ and $\psi_i(t)$ are orthogonal polynomials of suitable type in $[0,\ T]$, $T>0$.
To evaluate the unknowns $C_{i, j}$ we consider a finite number of
collocation points over the interval $[-1, 1]$.
For collocation
we have used the Tchebycheff nodes $
  x_j=\cos(\frac{\pi j}{N}), ~j=0,1,2,\cdots, M. $
Now collocating the differential operators at the Tchebycheff nodes $x = x_k$ yield
%
\[
  \frac{\partial^l}{\partial{t^l}} u_{h}(x_k,t)=\sum_{i=0}^{N} \sum_{j=0}^{M} C_{i j}\psi_i^{(l)}(t)\phi_{j} (x_k),\ k = 0, 1, \cdots, M,\  l=1,\ 2,
\]
and
\[
  \frac{\partial^l}{\partial{x^l}} u_{h}(x_k,t)=\sum_{i=0}^N\sum_{j=0}^{M} C_{i,j} \psi_i(t)\phi_{j}^{(l)} (x_k),\ k = 0, 1, \cdots, M,\  l=1,\ 2.
\]
Introducing time dependent parameters
\begin{equation}\label{timeparameters_01}
  a_j(t)= \sum_{i=0}^{N} C_{i j}\psi_i(t), \forall \ j=0,\ 1,\ \cdots, M,
\end{equation}
substituting  $\frac{\partial{ u_{h}(x_k,t)}}{\partial{t}}$, $\frac{\partial{ u_{h}(x_k,t)}}{\partial{x}}$, and
$\frac{\partial^2{ u_{h}(x_k,t)}}{\partial{x}^2}$ in \eqref{appr:f}, yields
the time  dependent system of linear differential equations
\begin{equation}\label{appr02:f}
  \left(\beta
 \frac{\partial \underline a(t)}{\partial t}+ \alpha\frac{\partial^2 \underline a(t)}{\partial t^2} \right)=  \mathbb{M}^{-1}A \underline a(t) +  \underline f(t),
\end{equation}
where $\underline a(t)= (a_0(t), \cdots, a_M(t))^T $ and $
\mathbb{M}_{i, j}= \phi_j(x_i), ~i, j=0,1,2,\cdots, M,~~$
$$
A_{i, j}  = \varepsilon \phi_j''(x_i) +  \mu \phi_j'(x_i) + \nu \phi_j(x_i),\ \ \mbox{$i, j=0,1,2,\cdots, M$}.
$$
Considering $
\Phi(x) = (\phi_0(x)\ \phi_1(x)\ \cdots\ \phi_M(x))', $
we write the boundary conditions as $
\underline a(t) \Phi(-1) =0, \ \text{and}\ \underline a(t) \Phi(1) =0. $
Since function values at both the boundary nodes are known (impose BCs $u(\pm 1)=0$),
ignoring elements that corresponds to $x_0$ and $x_M$ yields $\mathbb{M}$ and $A$ to the matrices of size $M-1\times M-1$.
Using \eqref{timeparameters_01} we rewrite
\[
  u_h(x_k, t_0)= \sum_{j=0}^N a_j(t_0) \phi_j(x_k)=u_0(x_k), \ \forall \  x_k.
\]
Therefore, \eqref{appr02:f} constitutes system of differential equations with initial functions
$ \underline a(t_0) = \mathbb{M}^{-1} G_0,~~\underline a'(t_0) = \mathbb{M}^{-1} G_1,
$ where $G_0 = \left[ g_0(x_1)\ g_0(x_2)\ \cdots g_0(x_{M-1})\right]'$ is a $M-1$ column vector, and $G_1$ is a $M-1$ column vector as $ G_1 = \left[ g_1(x_1)\ g_1(x_2)\ \cdots g_1(x_{M-1})\right]'.
$.

\section{Computation of parameters $C_{i, j}$}\label{compt}

Here we aim to compute $C_{i, j}$ values for the  approximation \eqref{appr1:f}. In this section we motivate ourselves to solve the time dependent system of linear differential equations \eqref{appr02:f} with given  initial conditions.
To that end we recall \eqref{appr02:f} with solutions as a truncated series of orthogonal polynomials given by
 \[
 a_i(t) = \sum_{j=0}^N C_{i,j} \psi_j(t), 0\le t\le T,
 \]
 which gives
  \[
 \left(
               \begin{array}{c}
                 a_0 (t) \\
                 a_1 (t) \\
                 . \\
                 . \\
                 . \\
                 a_N (t) \\
               \end{array}
             \right)
  = \left(
                                                                       \begin{array}{ccccccc}
                                                                          C_{00} & C_{01} & C_{02} & . & . & . & C_{0N} \\
                                                                         C_{10} & C_{11} & C_{12} & . & . & . & C_{1N} \\
                                                                         . & . & . & . & . & . & . \\
                                                                         . & . & . & . & . & . & . \\
                                                                         . & . & . & . & . & . & . \\
                                                                        C_{N0} & C_{N1} & C_{N2} & . & . & . & C_{NN} \\
                                                                       \end{array}
                                                                     \right)
 \left(
               \begin{array}{c}
                 \psi_0 (t) \\
                 \psi_1 (t) \\
                 . \\
                 . \\
                 . \\
                 \psi_N (t) \\
               \end{array}
             \right), \]
can be expressed as $ a_i(t) =  C_i  \psi(t) $, where $ \psi(t)  = \left[\psi_0(t)\ \psi_1(t)\ \cdots\ \psi_N(t)\right]^T, $
and $ C_i = [C_{i, 0}\ C_{i, 1}\ C_{i, 2} \cdots C_{i, N}]~;~i=0, 1,\ 2, \cdots, N$.
for all  $i=0, 1,\ 2, \cdots, N$, and $\psi_j(t)$ are orthogonal polynomials in $[0\ T]$.
Following \cite{Msezer} we write the derivatives of the unknown function $a_i(t)$ as
$$
   a_i^{(k)}(t) = \sum_{j=0}^N C_{i,j}^{(k)} \psi_j(t),\ k=1,\ 2,\ \cdots,
$$
and we rewrite the above expression as $a_i^{(k)}(t) = \psi(t) C_i^{(k)}.$ Also $C_i^{(k)}$ can be computed as
$ C_i^{(k)}  = 2^k \mathcal{M}^k C_i.$
Accordingly~\cite{Msezer},\\
\[ \text{when} ~ N ~\text{is}~ \text{odd},~~~
\mathcal{M}  = \left(
                \begin{array}{ccccccc}
                 0 & \frac{1}{2} & 0 & \frac{3}{2}& 0 &\cdots & \frac{N}{2}\\
                 0 & 0 & 2& 0& 4&0 \cdots & 0\\
                 0 & 0 & 0& 3 & 0 & 5\cdots & N\\
                 \vdots&\vdots&\vdots&\vdots &&&\vdots\\
                 0& 0 & 0 & & \cdots && N\\
                 0 & 0 & 0 & 0 & 0 \cdots&0 & 0
                \end{array}
                \right),
 \]
 and
 \[ \text{when}~ N ~\text{is}~ \text{even},~~~
\mathcal{M}  = \left(
                \begin{array}{ccccccc}
                 0 & \frac{1}{2} & 0 & \frac{3}{2}& 0 &\cdots & 0\\
                 0 & 0 & 2& 0& 4&0 \cdots & N\\
                 0 & 0 & 0& 3 & 0 & 5\cdots & 0\\
                 0& 0 & 0 & & \cdots && N\\
                 0 & 0 & 0 & 0 & 0 \cdots&0 & 0
                \end{array}
                \right).
 \]

 Thus the time dependent parameter and the derivatives can be written as
 \[
   a_i^{(k)}(t) = 2^k \psi(t)\mathcal{M}^k C_i,
 \]
and thus
 \[
 \underline a^{(k)}(t) = 2^k {  \Psi(t)} \mathcal{\mathbf{M}}^k C,\ k=0,\ 1,\ 2, \cdots,
 \]
where
 \[
   \Psi(t)  = \textbf{diag}( \psi(t), \psi(t), .~ .~ .~, \psi(t) )_{(M-1)\times (M-1)},
 \]

 \[
   C = ( C_1, C_2, .~ .~ .~,  C_{M-1})'_{(M-1)\times 1},
 \]
 and
 \[
  \mathcal{\mathbf{M}}^k  = \textbf{diag}
   ( \mathcal{M}^k , \mathcal{M}^k, .~ .~ .~ \mathcal{M}^k)_{(M-1)\times (M-1)}.
  \]
 Using the above notations the time dependent system of ODEs \eqref{appr02:f} can be written as
 \begin{equation}\label{flappr02:f}
  \alpha U^{(2)}(t)+\beta U^{(1)}(t))-\mathbb{M}^{-1}A U^{(0)}(t) = \underline f(t),
 \end{equation}
 where
 \[
   U^{0}(t) = \Psi(t) C,\  U^{1}(t) = 2 \Psi(t)\mathcal{\mathbf{M}} C, \ \text{and}\ U^{2}(t) = 2^2 \Psi(t)\mathcal{\mathbf{M}}^2 C.
 \]

To solve the time dependent algebraic system we collocate \eqref{flappr02:f} at translated Tchebychef nodes.
To that end we consider Tchebychev quadrature points for get a full discrete system of equations.
Let $t_j$, $j=0,\ 1,\ 2,\ \cdots, \ N$ be the number of Tchebychev points  on $[0,\ T]$.
Using the prescribed grid points we compute the time dependent system as
 \begin{equation}\label{f:full_discrete_01}
  \alpha P_2 \underline U^{(2)}+\beta  P_1 \underline U^{(1)})+  P_0 \underline U^{(0)} = \underline F,
 \end{equation}
 where
 \[
  \underline U^{(i)},\ i=0,\ 1,\ 2, \ \mbox{are $(M-1)(N+1)$ unknown column vector,}
 \]
 and
 \[
  \underline F =  [\underline f(t_0) \ \underline f(t_1)\cdots \ \underline f(t_N)]', \ \mbox{a $(M-1)(N+1)$ vector,}
 \]
 \[
 P_2=\left(
     \begin{array}{ccccc}
    I_{(M-1)\times (M-1)} &&&&\\
    & I_{(M-1)\times (M-1)}&&&\\
     && \ddots && \\
    &&&& I_{(M-1)\times (M-1)}
    \end{array}
 \right),
 \]
 $\mbox{a $(M-1)(N+1) \times (M-1)(N+1)$ matrix}$, $P_1=P_2$, and
 \[
 P_0=\left(
    \begin{array}{ccccc}
    \mathbb{M}^{-1}A_{(M-1)\times (M-1)} &&&&\\
    & \mathbb{M}^{-1}A_{(M-1)\times (M-1)}&&&\\
    && \ddots&&\\
    &&&& \mathbb{M}^{-1}A_{(M-1)\times (M-1)}
    \end{array}
 \right),
 \]
 $\mbox{a $(M-1)(N+1) \times (M-1)(N+1)$ matrix}.$
 Thus \eqref{f:full_discrete_01} gives a full discrete system of equation of the form
 \[
 \mathcal{W}\underline C = \underline F.
 \]
 Now from initial condition we get
 $
   \Psi(t_0) \underline C =\underline G.
 $
 Thus replacing the last row of the matrix  $ \mathcal{W}$ by $\Psi(t_0)$, and  last $P$ elements of $\underline F$ by $\underline G$ we get a system of equations of the form
 \[
 \tilde{\mathcal{W}} \underline C = \underline {\tilde F},
 \]
 which can be solved
 by using any standard linear system solver.

\section{Numerical experiments and discussions}\label{exprmt}
The space time orthogonal polynomial approximation presented here is very efficient to solve a certain class of
partial differential equations. Here we present some numerical examples to demonstrate the working of our scheme.
We use $\phi_j(x)$ as  Lagrange polynomials and $\psi_i(t)$ as Tchebychef polynomials
at a set of clustered grids.  We use Tchebychef nodes to define $\phi_j(x)$ and $\psi_i(t)$.
We use MATLAB for implementation purpose. We demonstrate the scheme by the following examples.
%
\begin{example}
%
Consider the initial boundary value problem (IBVP)
\[
 u_t = \varepsilon u_{xx}, \ u(\pm 1) = 0, \ u(x, 0) = u_0(x).
\]
Here we apply the scheme  \eqref{f:full_discrete_01} to approximate the solution of the IBVP.
 Figure~\ref{test_heat_01}
   shows the result.
\begin{figure}[ht!]
     \begin{center}
     \includegraphics[width=0.45\textwidth,height=6cm]{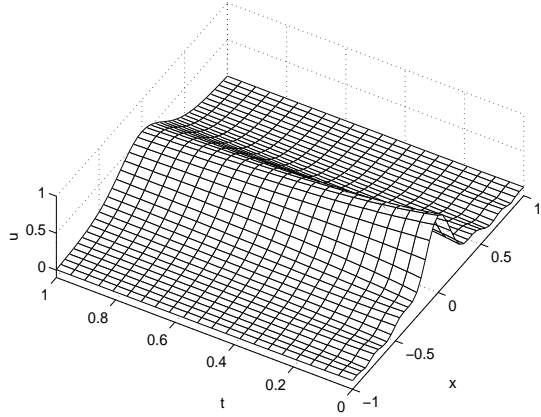}
     \end{center}
     \caption{Solution of the heat equation using \eqref{f:full_discrete_01}. Here we consider
      $u(x, 0) = e^{-5 |x|}$, $\varepsilon = 0.01$, $N=5$, $M=15$.
      }
     \label{test_heat_01}
\end{figure}
\end{example}
\begin{example}
We  consider  the IBVP
\[
u_t = \varepsilon u_{xx}+ \beta u_{x}+ \gamma u,  \ u(\pm 1) = 0, \ u(x, 0) = u_0(x)
\]
Here we apply the scheme  \eqref{f:full_discrete_01} to approximate the solution of the IBVP.
  Figure~\ref{test_advec_01}
   shows the results for various choices for $\beta$.
\begin{figure}[ht!]
     \begin{center}
       \includegraphics[width=0.49\textwidth,height=6cm]{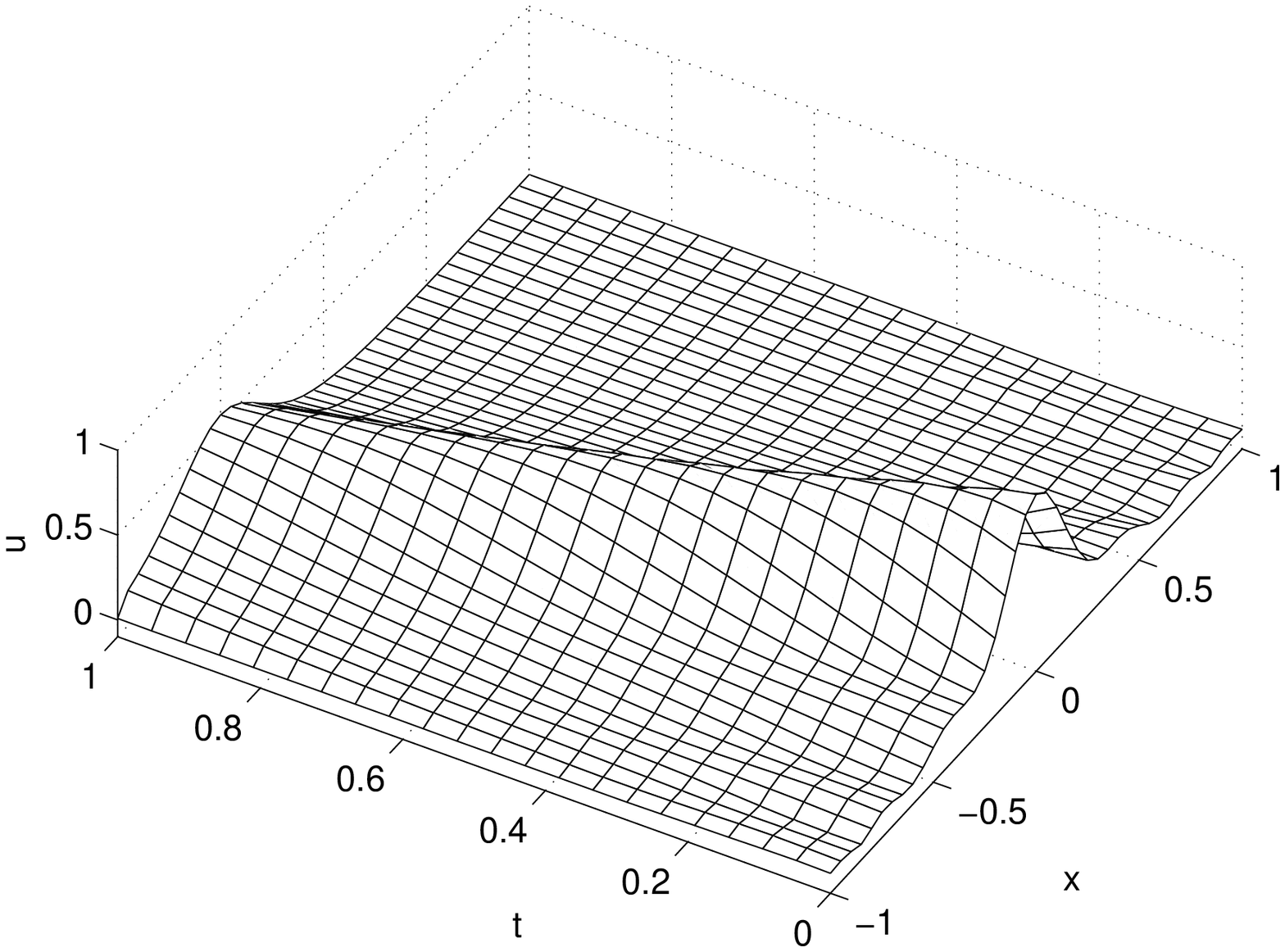}
       \includegraphics[width=0.49\textwidth,height=6cm]{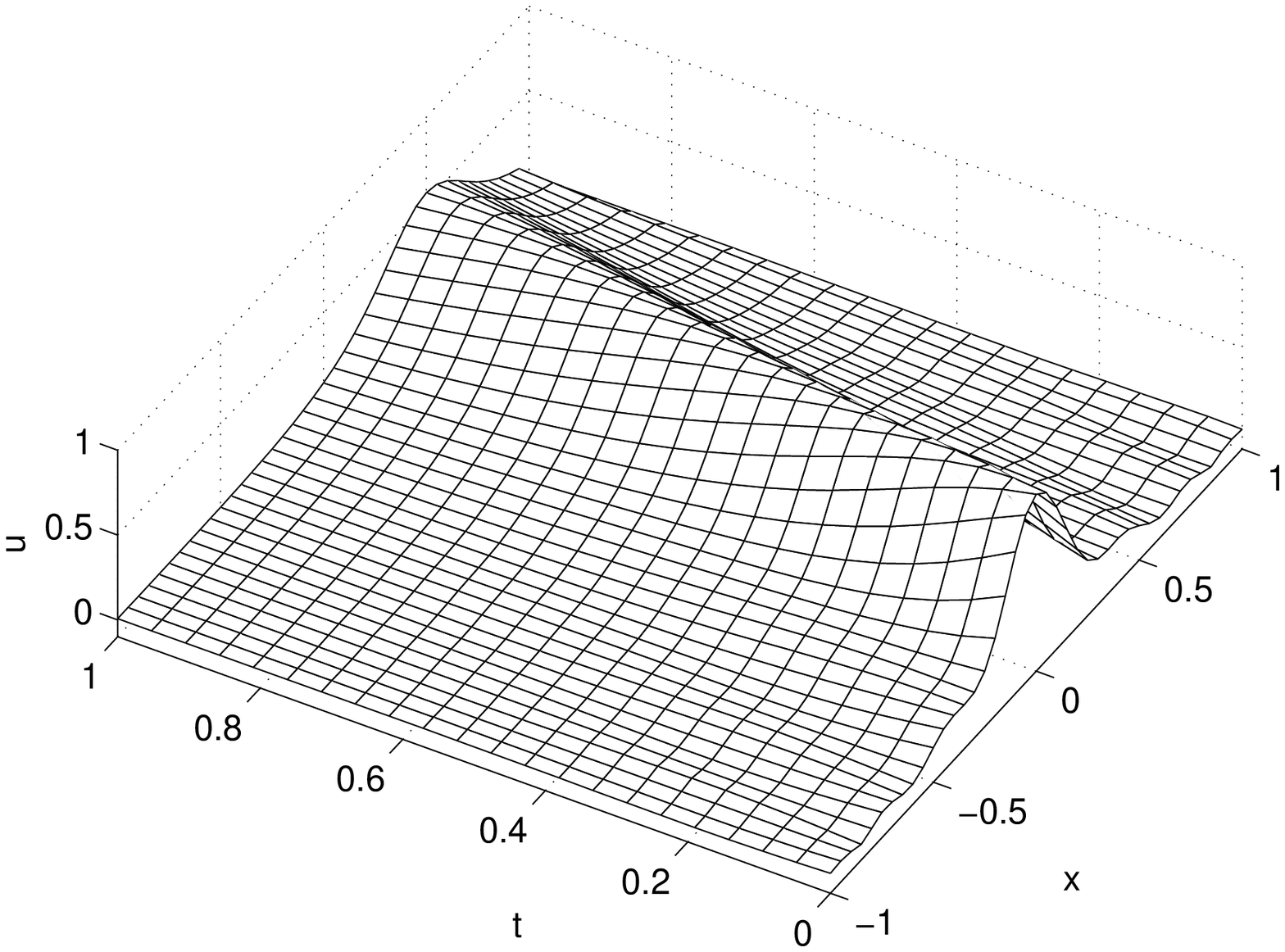}
     \end{center}
\caption{ Solution of the one way wave equation: pseudo spectral for space and  time integrations.
Here we consider
      $u(x, 0) = e^{-8\pi x^2}$,  $\varepsilon = 0.005$ and $\gamma=-0.1$, and $\beta =0.25$ (left),
      $\beta = -0.25$ (right) , $N=21$, $M=21$.
}
    \label{test_advec_01}
    \end{figure}

\end{example}
\begin{example}
We  consider  the IBVP
\[
  u_t = \varepsilon u_{x}, \ u(\pm 1) = 0, \ u(x, 0) = u_0(x)
\]
Here we apply the scheme  \eqref{f:full_discrete_01} to approximate the solution of the IBVP.
Figure~\ref{test_1wave_01} and
   Figure~\ref{test_1wave_02}
   show the results for various choices for $\varepsilon$.
\begin{figure}[ht!]
     \begin{center}
       \includegraphics[width=0.49\textwidth,height=6cm]{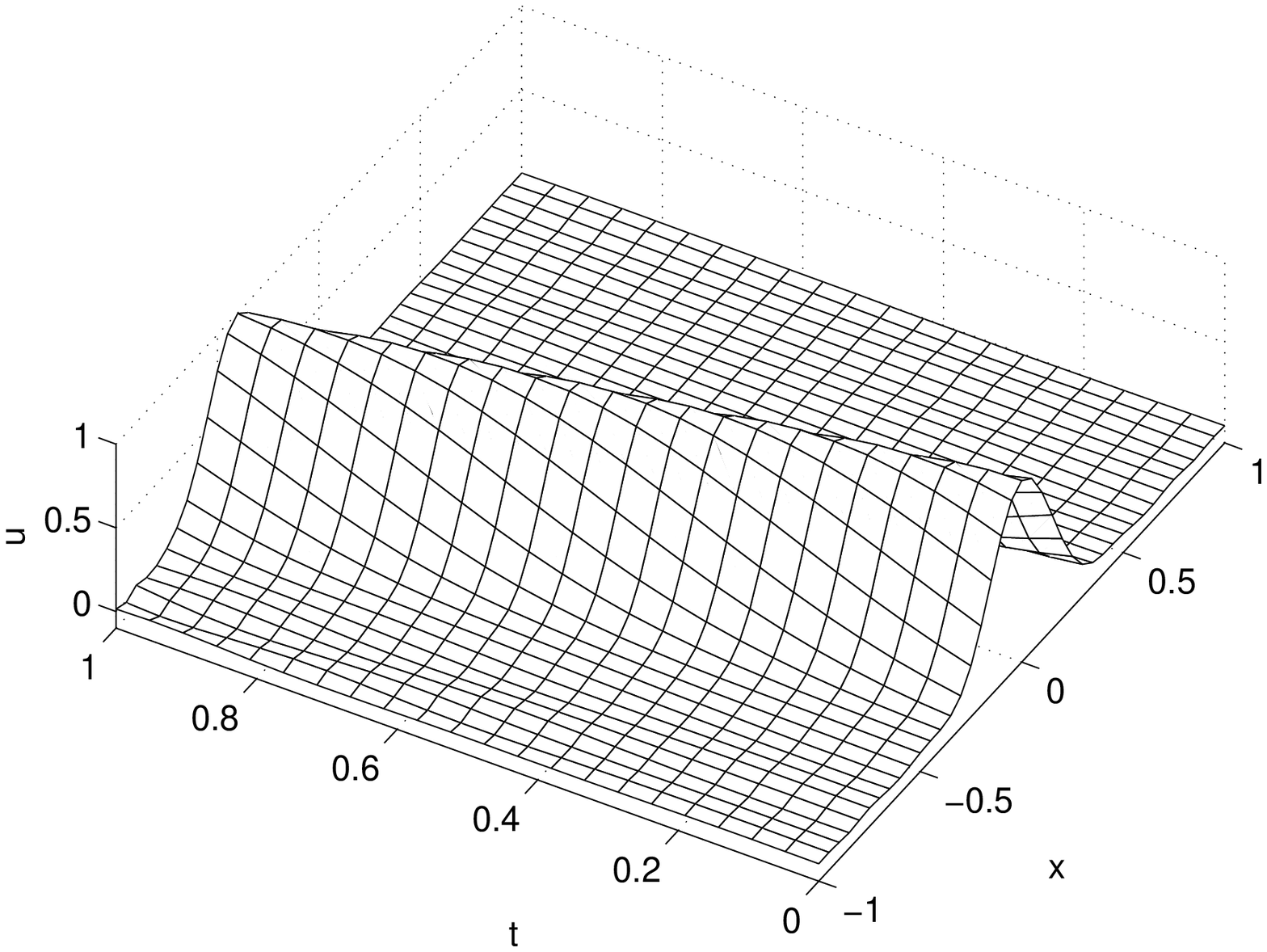}
       \includegraphics[width=0.49\textwidth,height=6cm]{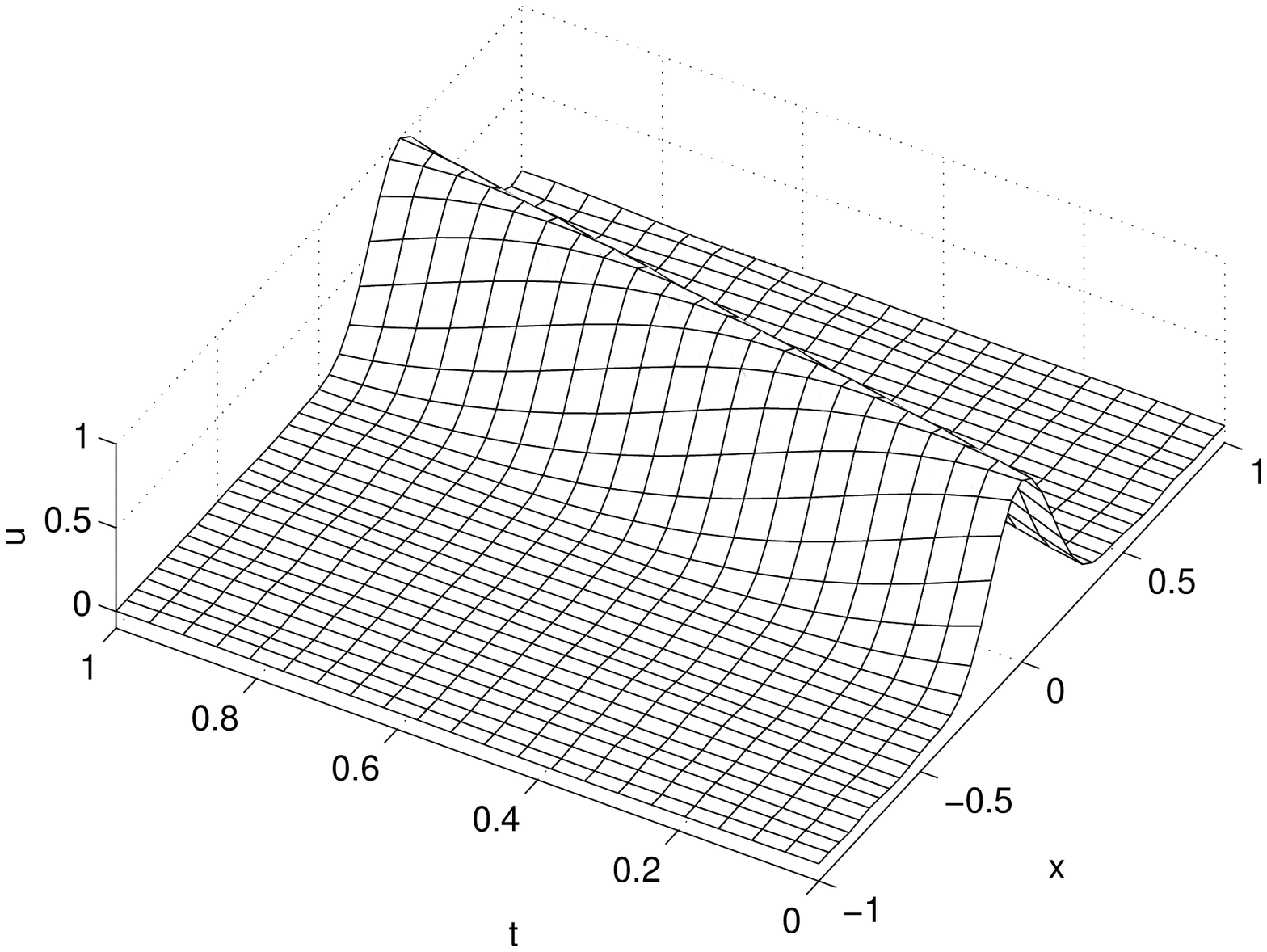}
     \end{center}
\caption{ Solution of the one way wave equation: pseudo spectral for space and  time integrations.
Here we consider
      $u(x, 0) = e^{-8\pi x^2}$, $\varepsilon = 0.2$(left),$\varepsilon = -0.2$(right) , $N=21$, $M=21$.
}
    \label{test_1wave_01}
    \end{figure}
\begin{figure}[ht!]
     \begin{center}
       \includegraphics[width=0.49\textwidth,height=6cm]{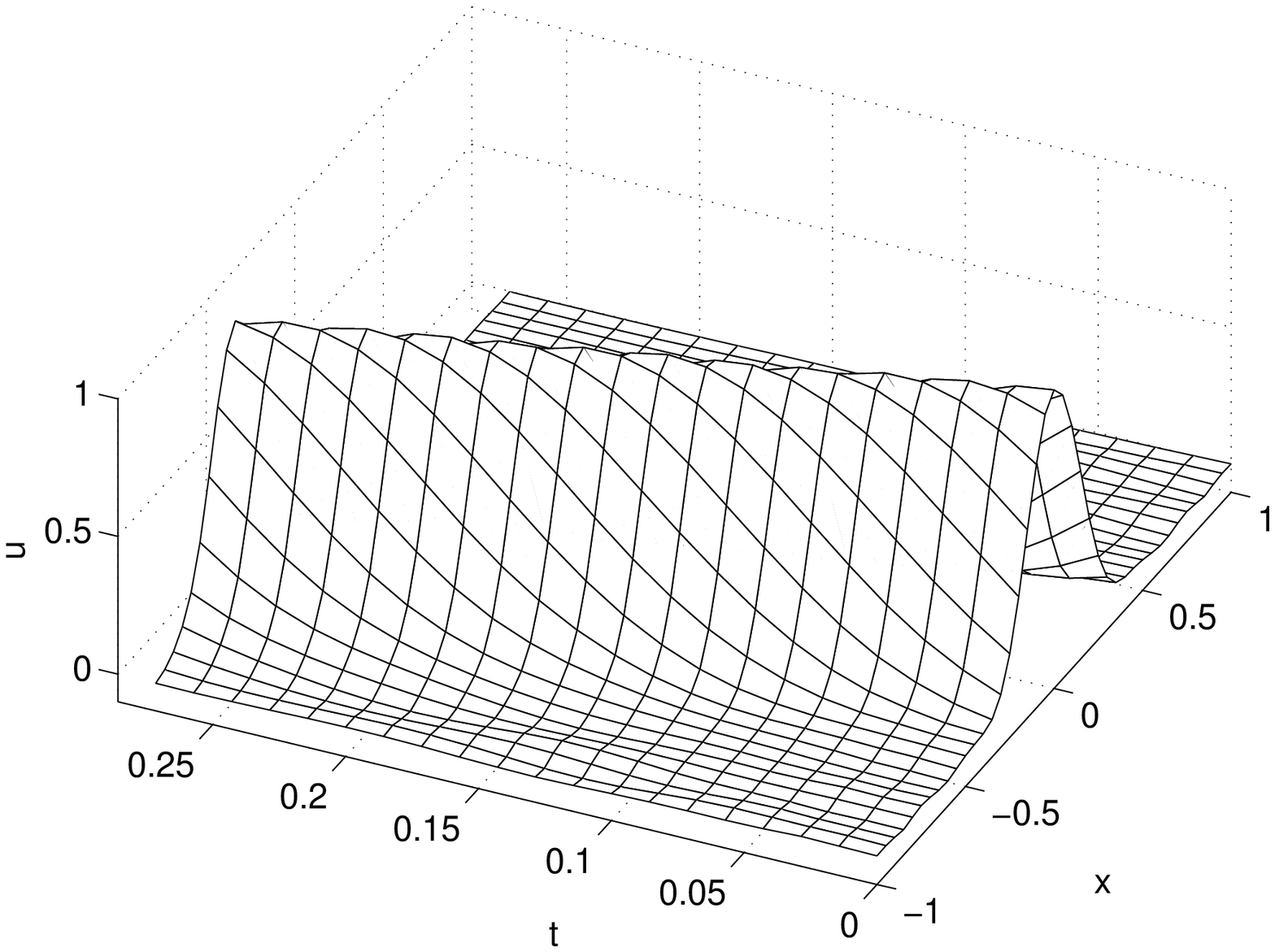}
       \includegraphics[width=0.49\textwidth,height=6cm]{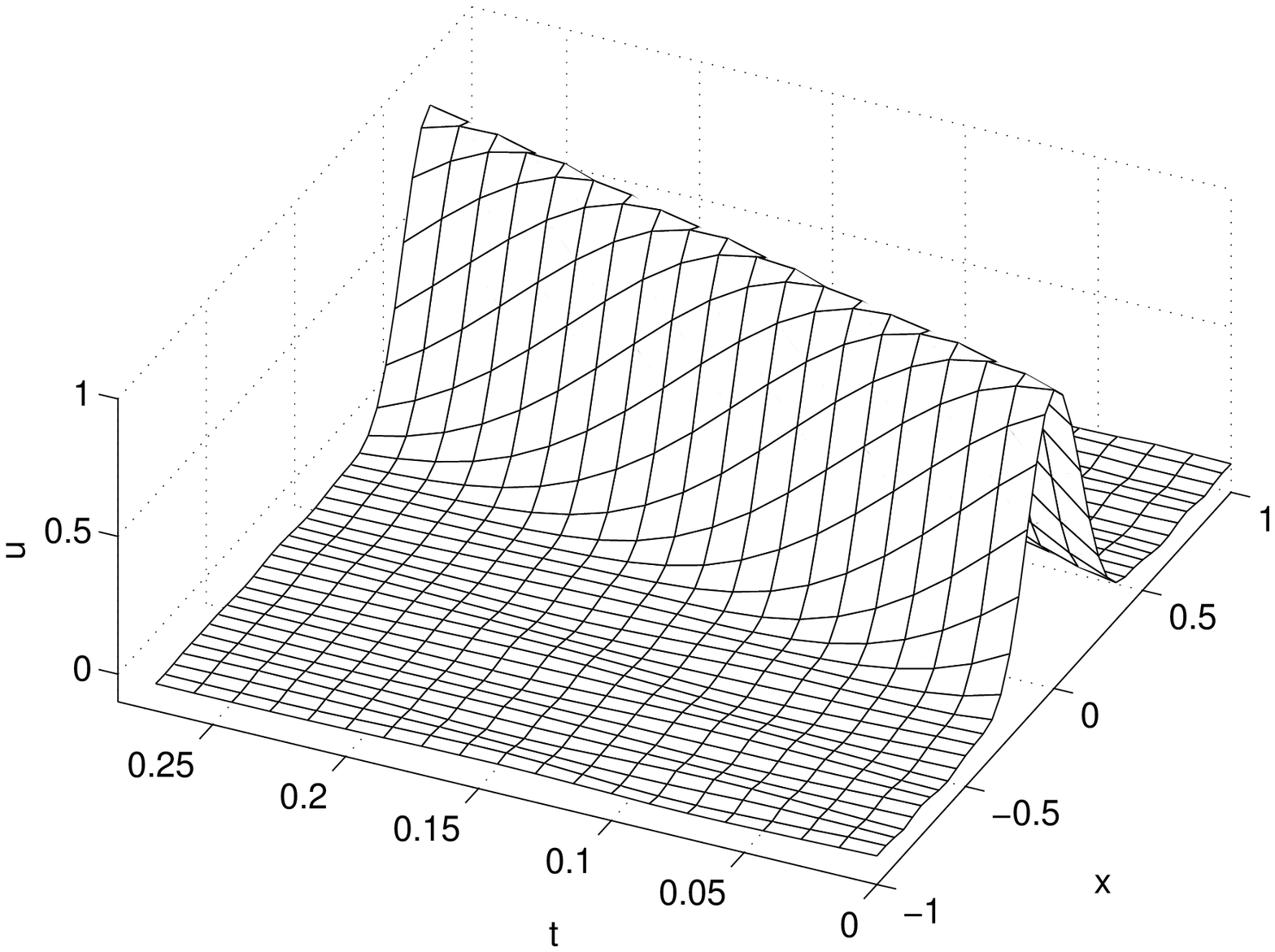}
       \includegraphics[width=0.49\textwidth,height=6cm]{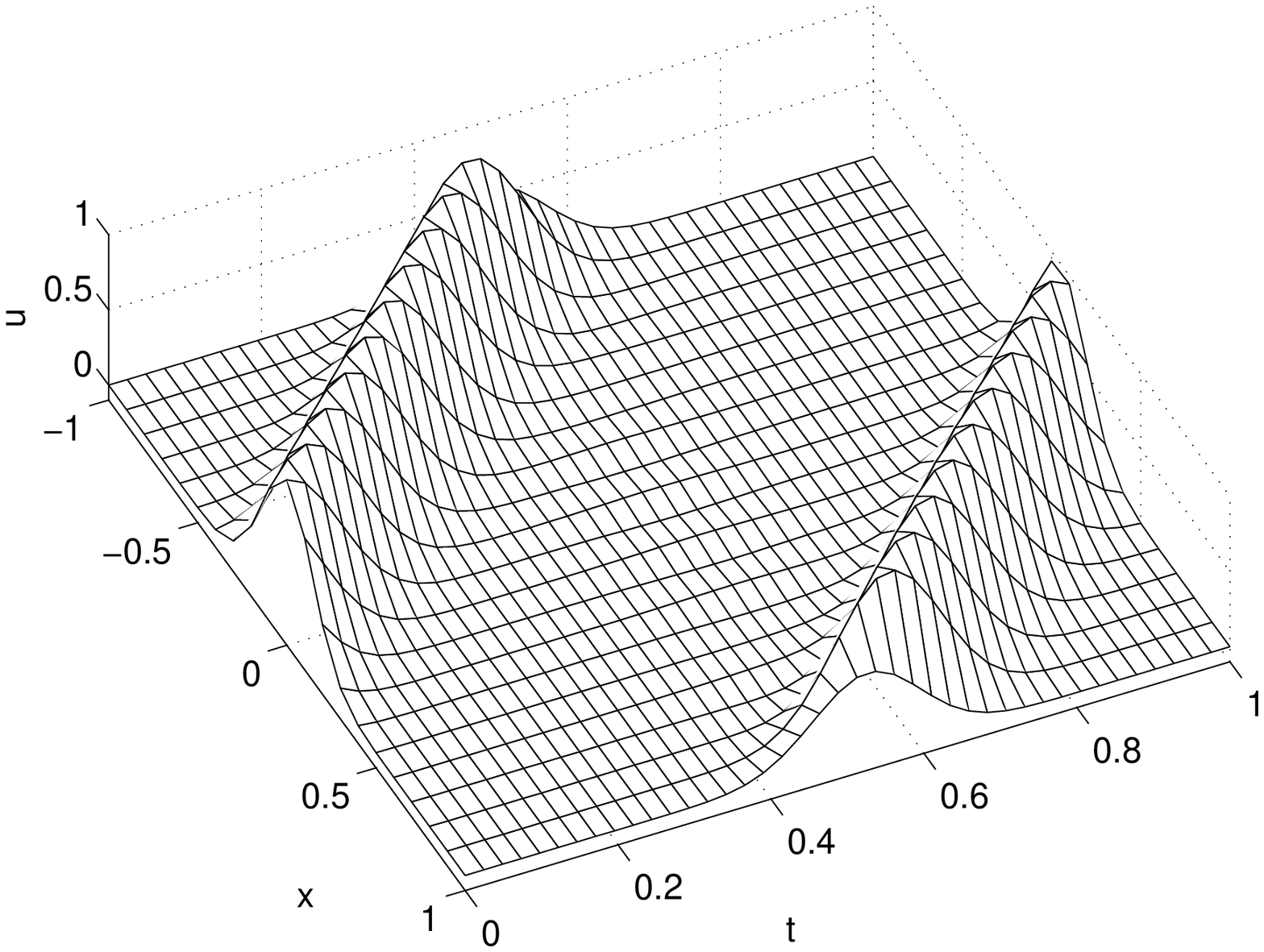}
       \includegraphics[width=0.49\textwidth,height=6cm]{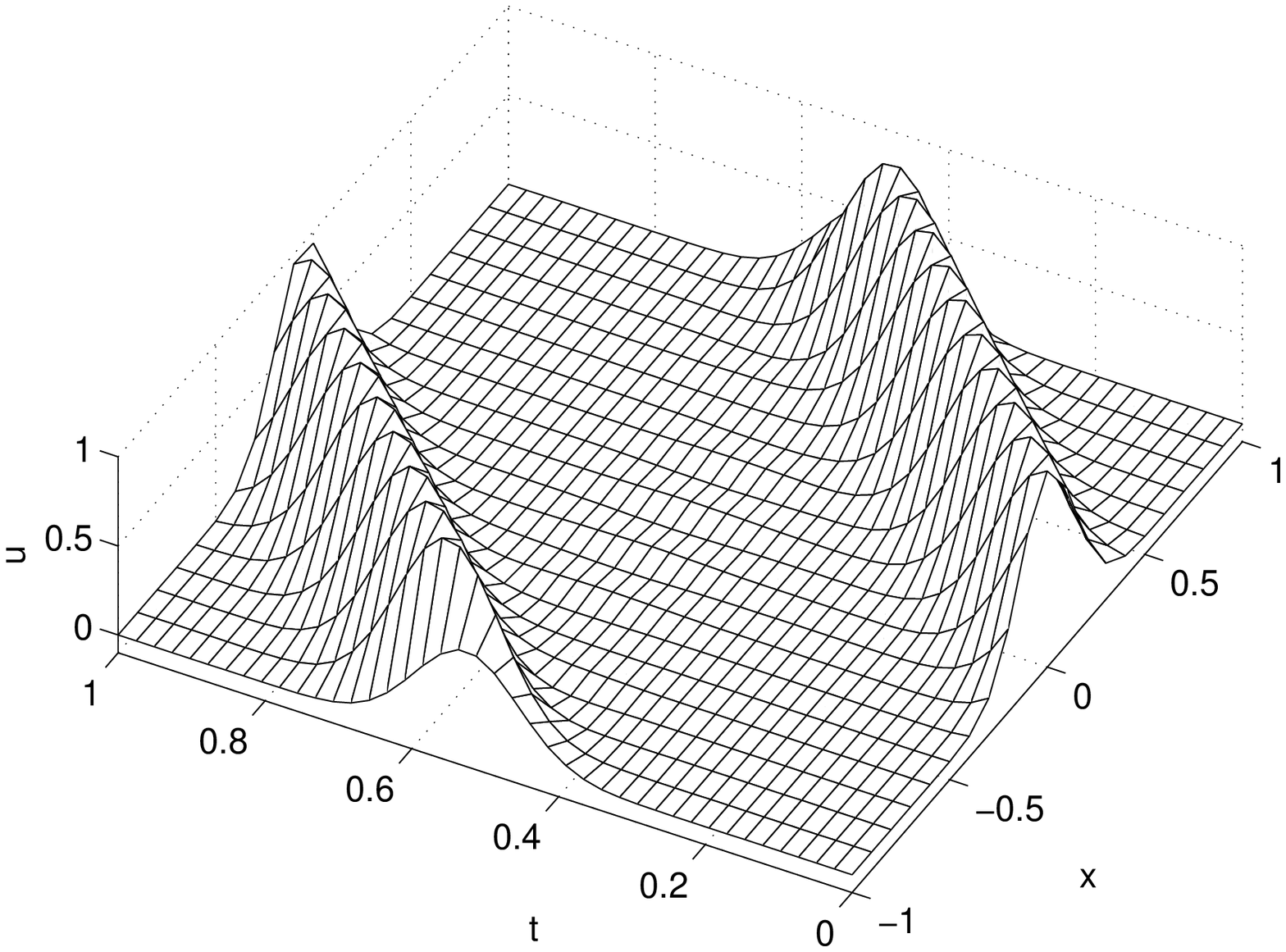}
     \end{center}
\caption{ Solution of the one way wave equation: pseudo spectral for space and  time integrations.
Here we consider $u(x, 0) = e^{-8\pi x^2}$, $\varepsilon = 1 $(left),$\varepsilon = -1$(right) ,
 $N=41$, $M=41$.  Upper figures show one fourth part of the solutions whereas the lower figures show solution for
 $ t \in
 [0,\ 1]$.
}
    \label{test_1wave_02}
    \end{figure}

\end{example}
There are many articles that approximate PDEs using  finite difference schemes and finite element schemes.
One needs a very fine grid points to get a reasonable solution using these schemes. That is to say one needs many grid points to approximate PDEs using FDS, or FES. Now from the examples presented above we notice that $N=21$ and $M=21$ are needed to approximate the solutions. Even a smaller choices of $N=5$, and $M=15$
can generate a reasonable approximate solution. As a result we need a very little computational time
to approximate the solutions. Whereas finite difference and finite element schemes generate a
huge matrix as a discrete partial differential operator. As a result these schemes are computationally expensive.
Thus space time orthogonal method can outperform most traditional schemes (finite difference scheme, finite element scheme, or spectral scheme for space followed by one/multistep scheme for temporal integration).

\section{Conclusion}\label{concl}
In this short communication, a computational scheme is developed for a class of PDEs using orthogonal polynomials. We carried out many numerical computations by taking Lagrange polynomial for space and Tchebycheb polynomials for time integration to validate this technique and demonstrate its capacity. We demonstrate that our scheme needs a very small number of grid points to compute the approximate solutions, whereas there are evidences~\cite{Atkinson, Trefethen} that FEM/FDM needs many grid points to compute solutions, thus they are computationally expensive. 
 Form this short study we conclude that the presented technique of polynomials for both space and time gives huge reduction of storage cost, and as a result it needs very small computational time. A future work is expected to work with bit more complex case with different boundary conditions and to extend the proposed computational scheme for some other multi-dimensional model problems.

\bibliography{ref}{}
\bibliographystyle{plain}
\end{document}